\documentclass{article}
\usepackage[utf8]{inputenc}

\usepackage{amsmath, amsthm, amsfonts, amssymb,latexsym, comment}
\newtheorem{lemma}{Lemma}[section]

\newtheorem{theorem}[lemma]{Theorem}
\newtheorem*{theorem-nonum}{Theorem}

\theoremstyle{remark}

\newcommand{\R}{\mathbb{R}}

\newcommand{\F}{\mathbb{F}}

\newcommand{\cK}{\mathcal{K}}

\usepackage[a4paper, total={6in, 8in}]{geometry}
\usepackage{color}
\usepackage[dvipsnames]{xcolor}

\newcommand{\new}[1]{{{\color{black}#1}}}
\newcommand{\cS}{{\cal S}}

\title{The Spherical Kakeya Problem in Finite Fields}
\author{Mehdi Makhul, Audie Warren and Arne Winterhof}

\begin{document}
\date{}
\maketitle

\begin{center}
Johann Radon Institute for Computational and Applied Mathematics,
Austrian Academy of Sciences, Altenberger Str. 69, 4040 Linz, Austria\\
E-mail: \{mehdi.makhul,audie.warren,arne.winterhof\}@oeaw.ac.at
\end{center}

\begin{abstract}
    We study subsets of the $n$-dimensional vector space over the finite field $\F_q$, for odd $q$, which contain either a sphere for each radius or a sphere for each first coordinate of the center. We call such sets radii spherical Kakeya sets and center spherical Kakeya sets, respectively.
    
    For $n\ge 4$ we prove a general lower bound on the size of any set containing $q-1$
    different spheres which applies to both kinds of spherical Kakeya sets. We provide constructions which meet the main
    terms of this lower bound.
   
     We also give a construction showing that we cannot get a lower bound of order of
     magnitude~$q^n$
    if we take lower dimensional objects such as circles in $\F_q^3$ instead of spheres, showing that there are significant differences to the line Kakeya problem.
    
    Finally, we study the case of dimension $n=1$ which is different and equivalent to the study of
    sum and difference sets that cover $\F_q$.
\end{abstract}

\let\thefootnote\relax\footnotetext{\emph{Mathematics Subject Classification}:
 52C10, 05B25, 11T99.}
\let\thefootnote\relax\footnotetext{\emph{Keywords and phrases.}
Kakeya problem, spheres, finite fields, diagonal equations, linear spaces.}

\section{Introduction}

A {\it (line-)Kakeya set} ${\cal K}\subset \F_q^n$ of $n$-dimensional vectors over the finite field $\F_q$ of $q$ elements is a set containing a line in each direction. It was shown
in \cite{dv} that every Kakeya set $\cal K$ satisfies $|{\cal K}|\ge c_nq^n$, where the implied constant $c_n$ depends only on the dimension $n$. Later research focused on the 
constant $c_n$, that is, on the one hand improved lower bounds \cite{dvkosasu} and on the other hand constructions of 'small' Kakeya sets \cite{kymuwa,ma,sasu}.

Several variants of Kakeya sets over finite fields have been studied as well, see for example \cite{elobta}.
In particular the paper \cite{wawi} deals with {\em conical Kakeya sets} over finite fields, that is, subsets of $\F_q^n$ containing either a parabola or a hyperbola in every direction (ellipses are not used since they do not have a direction). By 'directions' we usually mean points of the hyper-plane at infinity lying on an object.
This paper deals with spheres instead of lines. However, since spheres 
over finite fields have \emph{many} directions, 
roughly $q^{n-2}$ for $n \geq 3$, it is not desirable to use directions to define spherical Kakeya sets in finite fields. In analogy with the reals, we can define spherical Kakeya sets with reference to radii (see \cite{dehu,kowo,wo97} for real spherical Kakeya sets) or, say, the first coordinates of the centres of the spheres. 

Spheres over finite fields are well-studied objects, see \cite{ioko,iokosuthash,lini} and are defined as follows.
Throughout this paper we assume that $q$ is the power of an odd prime.
First we define the norm $\|\underline{x}\|$ of a vector in $\underline{x}=(x_1,\ldots,x_n)\in \F_q^n$ by
$$\|\underline{x}\|=x_1^2+\ldots+x_n^2.$$
In the finite field case this is more suitable than the square-root of the right hand side as used for the reals.
The {\em sphere} ${\cal S}_r(\underline{a})$ of radius $r\in \F_q^*$ and center $\underline{a}=(a_1,\ldots,a_n)\in \F_q^n$ is 
$${\cal S}_r(\underline{a})=\left\{\underline{x}\in \F_q^n : \|\underline{x}-\underline{a}\|=r\right\},$$
that is the set of solutions $\underline{x}=(x_1,\ldots,x_n)\in \F_q^n$ of the {\em quadratic diagonal equation}
$$(x_1-a_1)^2+\ldots+(x_n-a_n)^2=r.$$
Again in the finite field case it is more suitable to use $r$ instead of $r^2$ as in the real case.

Now a {\em radius spherical Kakeya set} in $\F_q^n$, $n\ge 2$, contains a sphere for each radius $r\in \F_q^*$ and a {\em (first coordinate of the) center spherical Kakeya set} in $\F_q^n$, $n\ge 2$,
contains a sphere for each first coordinate $a_1\in \F_q$ of the center.

For $n\ge 4$ we prove a general lower bound on sets ${\cal K}\subset\F_q^n$ which contain
$q-1$ different spheres which is also a lower bound on the size of spherical Kakeya sets.
We also provide a slightly different lower bound for $n=2,3$.

\begin{theorem}\label{dsphere}

Let $q$ be odd and $\cK \subset \F_q^n$ be a set containing at least $q-1$  distinct spheres for $n\ge 4$, or at least $(q-1)/2$ distinct
spheres for $n=2,3$. Then we have 
$$|\cK| \geq \begin{cases}\frac{1}{2}q^n + \frac{1}{2}q^{n-1} - q^{n-2} - \frac{1}{2}q^{\lfloor \frac{n-1}{2}\rfloor + 2} + \frac{1}{2} q^{\lfloor \frac{n-1}{2}\rfloor + 1},
\quad n \geq 4, \\ \frac{q^n-q^{n-2}}{4}, \quad n = 2,3. \end{cases}$$
\end{theorem}
In Section~\ref{lower} we prove Theorem~\ref{dsphere} by combining a well-known result on the number of solutions of quadratic diagonal equations with a simple counting argument. 

In Section~\ref{upper} we provide constructions of both radius spherical Kakeya sets and center spherical Kakeya sets which attain the main terms of this bound. In particular, we construct
a radius spherical Kakeya set of size
$$\frac{1}{2}q^n + \frac{1}{2}q^{n-1} -q^{n-2} + O\left(q^{n-3}\right)\quad \mbox{for }n\ge 8$$
and a center spherical Kakeya set of size
$$\frac{1}{2}q^n+\frac{1}{2}q^{n-1}+O\left(q^{n-2}\right)\quad \mbox{for }n\ge 5.$$
(We use the notation $X= O(Y)$ if $|X| \le cY$ for some absolute constant $c>0$.)

Now we introduce lower dimensional hyper-spheres, the motivation for which will be given in the next paragraph. 
Let ${\cal V}_{\underline{d}}=\{\underline{\new{x}}\in \F_q^n : \underline{d}\cdot \underline{x}=0\}$ be a linear subspace of $\F_q^n$ of dimension $n-1$ for some {\em direction} $\underline{d}\in \F_q^n\setminus\{\underline{0}\}$.
(We may assume that the first non-zero coordinate of $\underline{d}$ is $1$.)
Then the {\em hyper-sphere} ${\cal H}_{r}(\underline{a},\underline{d})$ in the hyper-plane $\underline{a}+{\cal V}_{\underline{d}}$ of radius $r\in \F_q^*$, direction $\underline{d}\in \F_q^n$ and center $\underline{a}\in \F_q^n$ is given by 
$${\cal H}_r(\underline{a},\underline{d})={\cal S}_r(\underline{a})\cap (\underline{a}+{\cal V}_{\underline{d}}).$$

In Section~\ref{hyper} we give a negative answer to the question of whether we could use lower-dimensional objects, for example circles in $\F_q^3$ instead of spheres, to get lower bounds of order of magnitude $q^n$. This question is motivated by the fact that the line Kakeya problem always deals with objects of dimension 1 (lines). However in our case, even hyper-spheres (which are of dimension $n-2$) are not enough to give asymptotic growth of order $q^n$. In particular, we show that in $\F_q^n$ there is a set of size $q^{n-1}+O(q^{n-2})$, $n\ge 3$, which contains a hyper-sphere for each center, direction and radius.

As in the real case \cite{dehu} our definition for spherical Kakeya sets in $\F_q^n$ can be adjusted for dimension $n=1$. A
circle ${\cal C}=\new{\{x\in \F_q: (x-a)^2=r^2 \}}$ in $\F_q$, for some radius $r\in \F_q^*$ and center $a\in \F_q$,  contains exactly two points $a\pm r$.
Note that here it is more suitable to use $r^2$ instead of $r$ (as for real circles).
A {\em radius circular Kakeya set in $\F_q$} contains a circle for each radius $r\in \F_q^*$, or equivalently we have
$$
{\cal K}-{\cal K}=\F_q,
$$
where
\begin{equation}\label{ominus}
{\cal K}-{\cal K}=\{x_1-x_2 : x_1,x_2\in {\cal K}\}.
\end{equation}
A {\em center circular Kakeya set in $\F_q$} contains a circle for each center $a\in \F_q$,
or equivalently we have
$${\cal K}\oplus{\cal K}=\F_q,$$
where
\begin{equation}\label{oplus}
{\cal K}\oplus {\cal K}=\{x_1+x_2: x_1,x_2\in {\cal K},x_1\ne x_2\}.
\end{equation}
In Section~\ref{dim1} we provide constructions of both radius circular and center circular Kakeya sets in $\F_q$ of optimal order of magnitude $O\left(q^{1/2}\right)$.

\section{Proof of Theorem~\ref{dsphere}}
\label{lower}

In this section we give a proof of Theorem~\ref{dsphere}
which is based on the following lemma.

\begin{lemma}\label{intersect}
The intersection of two different spheres $S_{r_1}(\underline{a})$
and $S_{r_2}(\underline{b})$, $(\underline{a},r_1)\ne (\underline{b},r_2)$, in $\F_q^n$, \new{where~$q$ is odd and~}$n\ge 2$, contains at most
$$q^{n-2}+q^{\lfloor (n-1)/2\rfloor}$$
points.
\end{lemma}
Proof. For $n\ge 1$, $a_1,\ldots,a_n\in \F_q^*$ and $r\in \F_q$ we recall that the number $N$ of solutions $(x_1,\ldots,x_n)\in \F_q^n$ 
to the quadratic diagonal equation
$$a_1x_1^2+\ldots+a_nx_n^2=r$$
satisfies
\begin{equation}\label{diag} |N-q^{n-1}|=  \left\{\begin{array}{cc} q^{\lfloor(n-1)/2\rfloor}, & r\ne 0,\\
 q^{\lfloor n/2\rfloor}-q^{\lceil (n-2)/2\rceil}, & r=0,\end{array}\right.
\end{equation}
see for example \cite[Theorem 10.5.1]{beevwi} or \cite[Theorems~6.26 and 6.27]{lini}.

For $n\ge 2$ we count the number of joint solutions $\underline{x}\in \F_q^n$ of the two equations
\begin{equation}\label{circ1} \|\underline{x}-\underline{a}\|=r_1
\end{equation}
and 
\begin{equation}\label{circ2} \|\underline{x}-\underline{b}\|=r_2.
\end{equation}
Subtracting $(\ref{circ2})$ from $(\ref{circ1})$ we get
\begin{equation}\label{lin}
2(\underline{b}-\underline{a})\cdot\underline{x}
=2(\underline{b'}-\underline{a'})\cdot\underline{x'}+2(b_n-a_n)x_n
=r_1-r_2-\|\underline{a}\|+\|\underline{b}\|,
\end{equation}
where $\underline{a}=(\underline{a'},a_n)$, $\underline{b}=(\underline{b'},b_n)$ and $\underline{x}=(\underline{x'},x_n)$ with $\underline{a'},\underline{b'},\underline{x'}\in \F_q^{n-1}$
and $a_n,b_n,x_n\in \F_q$.

If $\underline{a}=\underline{b}$ and thus $r_1\ne r_2$, then the two spheres are disjoint. Therefore we may assume
$\underline{a}\ne \underline{b}$. WLOG we may assume $a_n\ne b_n$.
\new{Then $x_n$ is of the form
$$x_n=\underline{u}\cdot \underline{x'}+c$$
by \eqref{lin}, where 
$$\underline{u}=(b_n-a_n)^{-1}(\underline{a'}-\underline{b'})$$
and 
$$c=(2(b_n-a_n))^{-1}(r_1-r_2-\|\underline{a}\|-\|\underline{b}\|).$$
Then we 
substitute $x_n$ in $(\ref{circ1})$ and get a quadratic form in at most $n-1$ variables,
$$\|\underline{x}-\underline{a}\|=\|\underline{x'}-\underline{a'}\|+(\underline{u}\cdot \underline{x'}+c-a_1)^2=r_1.$$
By \cite[Theorem~6.21]{lini} each quadratic form is equivalent to a diagonal equation, that is, it can be transformed into a diagonal equation by regular linear variable substitution. 
Hence, it}
has at most
$q^{n-2}+q^{\lfloor(n-1)/2\rfloor}$ solutions by $(\ref{diag})$ (applied with $n-1$ instead of $n$)
and the result follows.~\hfill $\Box$\\

We now prove Theorem \ref{dsphere}. Let ${\cal K}\subset \F_q^n$ contain at least $M$ different spheres ${\cal S}_1,\ldots,{\cal S}_M$.
By Lemma~\ref{intersect} each pair of spheres intersects in at most $q^{n-2}+q^{\lfloor (n-1)/2\rfloor}$ points, and each contains at
least $q^{n-1}-q^{\lfloor (n-1)/2\rfloor}$ points by $(\ref{diag})$.
Hence, 
$$\sum_{1\le i < j\le M} |S_i \cap S_j| \leq \left(q^{n-2}+q^{\lfloor (n-1)/2\rfloor}\right) \frac{M(M-1)}{2}$$
and we get
$$|{\cal K}|\ge \left|\bigcup_{i=1}^M{\cal S}_i\right|
\ge M\left(q^{n-1}-q^{\lfloor (n-1)/2\rfloor}\right)-\left(q^{n-2}+q^{\lfloor (n-1)/2\rfloor}\right)\frac{M(M-1)}{2}.$$
Choosing
$$M=\left\{\begin{array}{cc} (q-1)/2, & n=2 \mbox{ or }3,\\
                            q-1,      & n\ge 4,\end{array}\right.$$
we get
$$|{\cal K}|\ge  \frac{1}{2}q^n + \frac{1}{2}q^{n-1} - q^{n-2} - \frac{1}{2}q^{\lfloor \frac{n-1}{2}\rfloor + 2} + \frac{1}{2} q^{\lfloor \frac{n-1}{2}\rfloor + 1}
\quad \mbox{for }n \geq 4,$$ and 
$$|{\cal K}|\ge  \frac{q^n-q^{n-2}}{4}\quad \mbox{for }n = 2,3,$$
which completes the proof. \hfill $\Box$

\section{Constructions}
\label{upper}

In this section we give constructions of sets $\cK \subset \F_q^n$ containing either a sphere of every radius, or of~$q$ different first coordinates of the centres. In particular, for $n\geq8$, our construction for radii meets the constants in Theorem \ref{dsphere} up to and including the third term, and for $n\ge 5$, our construction for centers meets the first two constants.

\subsection{Spheres with different radii}

First we give a construction for different radii. 
For $r\in \F_q^*$ consider the sphere
$${\cal S}_r=\{\left(x,\underline{y}\right)\in \F_q^n : (x-r)^2+\|\underline{y}\|=r\}.$$
The union $\bigcup_{r \in \F_q^*} {\cal S}_r$ contains a sphere of every radius. We use the inclusion-exclusion principle to bound the size of this set. We firstly bound the intersection of two different spheres $\cS_r$ and $\cS_s$;
the intersection points are
$${\cal S}_r\cap {\cal S}_s=\left\{\left(\frac{r+s-1}{2},\underline{y}\right) : \|\underline{y}\|= rs - \left( \frac{r+s-1}{2}\right)^2\right\},\quad r\ne s,\quad  r,s \in \F_q^*.$$

$|{\cal S}_r\cap {\cal S}_s|$ is precisely the number of solutions $(y_1,...,y_{n-1})$ to the equation
$$y_1^2 + ... + y_{n-1}^2 =  rs - \left( \frac{r+s-1}{2}\right)^2.$$
Therefore, for each valid choice of $(r,s)$, we have $|\cS_r \cap \cS_s| = q^{n-2} + O\left(q^{\frac{n-1}{2}}\right)$ by $(\ref{diag})$. 
We can now explicitly find the sum of the size of intersections of any two spheres, as 
$$\sum_{\substack{r,s \in \F_q^* \\ r \neq s}}|\cS_r \cap \cS_s| = (q-1)(q-2)\left(q^{n-2} + O\left(q^{\frac{n-1}{2}}\right)\right) = q^n - 3q^{n-1} + 2q^{n-2} + O\left(q^{\frac{n+3}{2}}\right).$$
We can see via the $x$ coordinate $\frac{r+s-1}{2}$ that the intersection of any three distinct spheres $\cS_r$, $\cS_s$, and $\cS_t$ is empty. By the inclusion exclusion principle 
and $(\ref{diag})$
\begin{align*}\left|\bigcup_{r \in \F_q^*} {\cal S}_r\right| &= \sum_{r \in \F_q^*}|\cS_r| - \frac{1}{2} \sum_{\substack{r,s \in \F_q^* \\ r \neq s}}|\cS_r \cap \cS_s| \\
& =  \frac{1}{2}q^n + \frac{1}{2}q^{n-1} -q^{n-2} + O\left(q^{\frac{n+3}{2}}\right).
\end{align*}

\subsection{Spheres with different first coordinates of the centres}

For a fixed non-square $r\in \F_q^*$ consider the set 
$${\cal  Q}=\{\left(x,\underline{y}\right)\in \F_q\times \F_q^{n-1} : r-\|\underline{y}\| \mbox{ is a square in }\F_q\}.$$ 
The $q$ distinct spheres ${\cal S}_r(a)=\{\left(x,\underline{y}\right)\in \F_q\times \F_q^{n-1}: (x-a)^2+\|\underline{y}\|=r\}$, $a\in \F_q$,
are all subsets of~${\cal Q}$.
However, the size of ${\cal Q}$ is 
$$\frac{q^n+q^{n-1}}{2}+O\left(q^{n-2}\right)\quad \mbox{for }n\ge 5.$$ Indeed, by $(\ref{diag})$ each non-zero value of $\|\underline{y}\|$ is
attained 
$$q^{n-2}+O\left(q^{\lfloor n/2\rfloor -1}\right)$$ times, that is
$q^{n-2}+O(q^{n-4})$ for $n\ge 5$.
There are $(q+1)/2$ \new{(non-zero)} values $\|y\|$ such that $r-\|\underline{y}\|$ is a  square
(since $r$ is a non-square)
and $x$ can take any value in $\F_q$.

For $n=3,4$ we have
$$|{\cal Q}|=\frac{q^n}{2}+O\left(q^{n-1}\right).$$

\section{Hyper-spheres}
\label{hyper}

In this section we show for $n\ge 3$ that even if a set contains hyper-spheres for all directions, non-zero centers and radii, it may have only $q^{n-1}+O(q^{n-2})$ points.

We consider the union
$${\cal H}=\bigcup_{\underline{a}\in \F_q^n\setminus\{\underline{0}\}}{\cal H}_{-\|a\|}(\underline{a},\underline{a})$$
of the hyper-spheres 
$${\cal H}_{-\|a\|}(\underline{a},\underline{a})=\{\underline{x}\in \F_q^n : \|\underline{x} - \underline{a}\|+\|\underline{a}\|=\underline{a}\cdot (\underline{x}-\underline{a})=0\}, \quad \underline{a}\in \F_q^n\setminus \{\underline{0}\},$$ 
with center $\underline{a}$, direction $\underline{a}$ and
radius $-\|\underline{a}\|$ (which covers all radii since each 
element of $\F_q^*$ is sum of two squares). However,
each $\underline{x}\in {\cal H}_{-\|a\|}(\underline{a},\underline{a})$ satisfies
$$\|\underline{x}\|=\|\underline{x}-\underline{a}+\underline{a}\|=
\|\underline{x}-\underline{a}\|+\|\underline{a}\|+2\underline{a}\cdot(\underline{x}-\underline{a}) =0,$$
which has at most $q^{n-1}+q^{\lfloor n/2\rfloor}-q^{\lceil (n-2)/2\rceil}$ solutions by $(\ref{diag})$ which is an upper bound for $|{\cal H}|$.

\section{One-dimensional circular Kakeya sets}\label{dim1}

The definitions of circular Kakeya sets in dimension 1 are in fact equivalent to definitions concerning sum and difference sets. More precisely,
$\cK\subset \F_q$ is a {\it radius circular Kakeya set in $\F_q$} if and only if 
$${\cal K}-{\cal K} = \F_q$$
and
{\it a centre circular Kakeya set in $\F_q$} if and only if 
$${\cal K}\oplus{\cal K} = \F_q,$$
\new{where ${\cal K}-{\cal K}$ and ${\cal K}\oplus {\cal K}$ are defined by \eqref{ominus} and \eqref{oplus}.}

To see the first equivalence, let ${\cal K}\subset \F_q$ be a set that contains a circle of radius $r$ for each $r \in \F_q^*$. Therefore there exists $a\in \F_q$ such that $\left\{a+r, a-r\right\} \subset \cal{K}$. We get $a+r-(a-r)=2r \in \cal{K}-\cal{K}$. Therefore, since $2r$ covers all of $\F_q^*$, we have ${\cal K}-{\cal K}=\F_q$ ($0\in {\cal K}-{\cal K}$ trivially). Conversely, suppose that ${\cal K} \subset \F_q$ is a subset such that ${\cal K}-{\cal K}=\F_q$. Then for each $r \in \F_q$, there exist $x_1,x_2 \in \cal{K}$, such that $x_1-x_2=2r$. By taking $a=(x_1+x_2)/2$ we see $x_1=a+r$ and $x_2=a-r$ and that the circle $\left\{a+r, a-r\right\}$ is in  $\cal{K}$. 

For the second equivalence, let ${\cal K} \subset \F_q$ be a set containing a circle for any center $a$. Then for all $a \in \F_q$, there exists $r\in \F_q^{*}$ such that $\left\{a-r,a+r \right\} \subset \cal{K}$. Then we have $(a-r)+(a+r)=2a$, and therefore ${\cal K}\oplus{\cal K}=\F_q$. Conversely, let $\cal{K}$ be a subset of $\F_q^*$ such that ${\cal K}\oplus{\cal K}=\F_q$. Fix $a\in \F_q$. Since ${\cal K}-{\cal K}=\F_q$, there exist $x_1, x_2 \in \cal{K}$, $x_1\ne x_2$, such that $x_1+x_2=2a$. Taking $r=(x_1-x_2)/2$ we can write $x_1=a + r$ and $x_2=a-r$, so that a circle of centre $a$ is in $\cK$. 

Since $|{\cal K}-{\cal K}|\le |{\cal K}|^2$, each radius circular Kakeya set in $\F_q$ has size at least $\lceil q^{1/2}\rceil$, and since $|{\cal K}\oplus{\cal K}|<|{\cal K}|^2/2$ the  size of any center circular Kakeya set ${\cal K}$  of $\F_q$ is at least $|{\cal K}|\ge \lceil \sqrt{2q}\rceil$.
\new{(Keep the condition $x_1\ne x_2$ in \eqref{oplus} in mind.)}
In this section we will give constructions of radius circular and center circular Kakeya sets ${\cal K}$ in $\F_q$ with $|{\cal K}|$
of optimal order of magnitude $O(q^{1/2})$.

For a prime $p>2$ it is easy to find circular Kakeya sets in $\F_p$ of size $2\lfloor\sqrt{p}\rfloor+1$, 
 \begin{equation}\label{simpleconstruction}{\cal K}={\cal K}_p=\{0,1,2,...,\lfloor\sqrt{p}\rfloor \} \cup -{\cal K}_0,\end{equation}
 where 
 $${\cal K}_0=\left\{\lceil\sqrt{p}\rceil, 2\lceil\sqrt{p}\rceil,...,\lfloor\sqrt{p}\rfloor\lceil\sqrt{p}\rceil\right\}.$$
It is clear that ${\cal K}-{\cal K} = \F_p.$ 
Substituting $-{\cal K}_0$ by ${\cal K}_0$ in $(\ref{simpleconstruction})$ we get
${\cal K}\oplus{\cal K}=\F_p$.

If $q=r^2$ is a square and $\alpha$ is a defining element of $\F_q$ over $\F_r$, that is, $\F_q=\F_r(\alpha)$, then we can choose 
$${\cal K}=\F_r\cup \alpha \F_r$$
of size $|{\cal K}|=2q^{1/2}-1$ to get both
${\cal K}-{\cal K}=\F_q$ and ${\cal K}\oplus{\cal K}=\F_q$.

If $q=p^{2m+1}$ with a prime $p$ and $\F_q=\F_p(\beta)$ with a defining element $\beta$ of $\F_q$ over $\F_p$, then we first choose the construction ${\cal K}_p$ from $(\ref{simpleconstruction})$ and then take
$${\cal K}={\cal K}_1\cup {\cal K}_2,$$
where 
$${\cal K}_1=\{a_0+a_1\beta+\ldots+a_m\beta^m: a_0\in {\cal K}_p, a_1,\ldots,a_m\in \F_p\}$$
and
$${\cal K}_2=\{a_0+a_1\beta^{m+1}+\ldots+a_m\beta^{2m}: a_0\in {\cal K}_p,a_1,\ldots,a_m\in \F_p\}.$$
It is easy to check that 
${\cal K}-{\cal K}=\F_q$ and 
$$|{\cal K}|=(2p^m-1)|{\cal K}_p|<4q^{1/2}+2(q/p)^{1/2}.$$
Again substituting $-{\cal K}_0$ by ${\cal K}_0$ in $(\ref{simpleconstruction})$ we get
${\cal K}\oplus{\cal K}=\F_q$.

Combining all the cases we can formulate a general result.

\begin{theorem}
 For a fixed power $q$ of an odd prime let ${\cal K}\subset \F_q$ be either a radius circular or a center circular Kakeya set in $\F_q$ of minimal size. Then we have
 $$q^{1/2}\le |{\cal K}|<6q^{1/2}.$$
\end{theorem}
The constants can be certainly improved using, for example, ideas from \cite{Fried1988,jish}. However, we did not calculate these improved constants for the readability of this paper, and
since even the improved upper bounds would not be optimal.

\section*{Acknowledgment}
The authors are supported by the Austrian Science Fund FWF Project P~30405-N32.
We would like to thank Oliver Roche-Newton for pointing out this problem.

\new{We wish to thank the anonymous referees for their careful study of our paper and their very useful comments. }


\begin{thebibliography}{99}
   
  \bibitem{beevwi} {\sc B. C. Berndt, R. J. Evans and K. S. Williams} {\em Gauss and Jacobi sums.}
Canadian Mathematical Society Series of Monographs and Advanced Texts. A Wiley-Interscience Publication. John Wiley \& Sons, Inc., New York, 1998.
  
  
  \bibitem{dehu}  {\sc Y. Deng, C. Hu, S. Long, T. Tang, J. Thuswaldner and L. Xi.} {\em On a variant of the Kakeya problem in $\R$.} Arch. Math. (Basel) 101 (2013), no. 4, 395--400.
  

\bibitem{dv} {\sc Z. Dvir.} {\em On the size of Kakeya sets in finite fields.} J. Amer. Math. Soc. 22 (2009), no. 4, 1093--1097.

\bibitem{dvkosasu} {\sc Z. Dvir, S. Kopparty, S. Saraf and M. Sudan.} {\em Extensions to the method of multiplicities, with applications to Kakeya sets and mergers.} SIAM J. Comput. 42 (2013), no. 6, 2305--2328.

\bibitem{elobta} {\sc J. S. Ellenberg, R. Oberlin and T. Tao.}
{\em The Kakeya set and maximal conjectures for algebraic varieties over finite fields. }
Mathematika 56 (2010), no. 1, 1--25. 



\bibitem{Fried1988}{\sc K. Fried.} {\em Rare bases for finite intervals of integers.} Acta Sci. Math. (Szeged) 52 (1988), no. 3--4, 303--305.


 \bibitem{ioko} {\sc A. Iosevich and D. Koh.}
{\em Extension theorems for spheres in the finite field setting. }
Forum Math. 22 (2010), no. 3, 457–-483.

\bibitem{iokosuthash} {\sc A. Iosevich, D. Koh, S. Lee, T. Pham and C. Shen.} {\em On restriction estimates for spheres in finite fields.} Preprint 2019.

 \bibitem{jish} {\sc X. Jia and J. Shen.} {\em Extremal bases for finite cyclic groups.} SIAM J. Discrete Math. 31 (2017), no. 2, 796--804.



\bibitem{kowo} {\sc L. Kolasa and T. Wolff.}
{\em On some variants of the Kakeya problem. }
Pacific J. Math. 190 (1999), no. 1, 111--154. 

\bibitem{kymuwa} {\sc G. Kyureghyan, P. M\"uller and Q. Wang.} {\em On the size of Kakeya sets in finite vector spaces.} Electron. J. Combin. 20 (2013), no. 3, Paper 36, 10 pp.

\bibitem{lini}  {\sc R. Lidl and H. Niederreiter.} {\em Finite fields.}  Second edition. Encyclopedia of Mathematics and its Applications, 20. Cambridge University Press, Cambridge, 1997.

\bibitem{ma} {\sc A. Maschietti. }{\em Kakeya sets in finite affine spaces.} J. Combin. Theory Ser. A 118 (2011), no. 1, 228--230.


\bibitem{sasu}  {\sc S. Saraf and M. Sudan.} {\em An improved lower bound on the size of Kakeya sets over finite fields.} Anal. PDE 1 (2008), no. 3, 375--379.


\bibitem{wawi} {\sc A. Warren and A. Winterhof.} {\em Conical Kakeya and Nikodym sets in finite fields.} Finite Fields Appl. 59 (2019), 185--198.


 \bibitem{wo97} {\sc T. Wolff.} {\em A Kakeya-type problem for circles.} Amer. J. Math. 119 (1997), no. 5, 985–-1026.
\end{thebibliography}
\end{document}